\documentclass[12pt]{article}
\usepackage{amsmath,amsthm,amssymb,amscd}
\usepackage{latexsym}
\newtheorem{thm}{Theorem}[section]
 
 \newtheorem{lem}[thm]{Lemma}
 
 \theoremstyle{definition}
 
 \theoremstyle{remark}

 \numberwithin{equation}{section}

\title
{The transverse Chern-Ricci flow}

\author{ Hong Huang}
\date{}
\begin{document}
\maketitle
\begin{abstract}
 We introduce transverse Chern-Ricci flow for transversely Hermitian foliations,  which is analogous to the Chern-Ricci flow. We show that when $\mathcal{F}$ is homologically orientable and the basic first  Bott-Chern class is zero, starting at any transversely Hermitian metric the flow exists for all time and as $t\rightarrow \infty$  converges smoothly to a transversely Hermitian metric $\omega_\infty$ with the transverse Chern-Ricci form $\rho^T(\omega_\infty)=0$. We also characterize  the maximal existence time of the flow in the general case. These are foliated version of  results of Gill and Tosatti-Weinkove, and also extend recent work of Bedulli-He-Vezzoni.


\end{abstract}
\maketitle


\section {Introduction}

In 2000, inspired by  Hamilton's  seminal work  on Ricci flow [H2],   Lovri$\acute{c}$, Min-Oo  and Ruh  [LMR] introduced the transverse Ricci flow

\begin{equation}
\frac{\partial  g^T}{\partial t}=-Ric^T, \hspace*{8mm}
g^T(0)=g_0^T
\end{equation}
 for a Riemannian foliation on a compact manifold $(M,g_0)$,
 and showed the short time existence and uniqueness following Hamilton's original approach [H2] via Nash-Moser inverse function theorem.  Here, $g^T$ and $Ric^T$ are the transverse metric and Ricci tensor respectively, viewed as sections of $S^2(Q^*)$, where $Q$ is the normal bundle to the foliation.

In the special case when $(M,g_0)$ is a compact Sasaki manifold with  the basic first Chern class $c_1^B(M)$  a multiple of the basic class of the transverse K$\ddot{a}$hler form, say $c_1^B(M)=\kappa [\frac{1}{2}d\eta_0]_B$, where $\kappa=1, 0$ or $-1$, Smoczyk, Wang and Zhang [SWZ] showed that the flow
\begin{equation}
\frac{\partial  g^T(t)}{\partial t}=-(Ric^T(t)-\kappa g^T(t)), \hspace*{8mm}
g^T(0)=g_0^T
\end{equation}
has a long time solution, and one can get Sasaki metrics $g(t)$ which induce the transverse metrics $g^T(t)$. (1.2) is called the Sasaki-Ricci flow.  In [SWZ] it is also showed that when $\kappa= 0$ or $-1$ the flow (1.2) converges to an $\eta$-Einstein metric. The Sasaki-Ricci flow  was further studied by Collins, He and some other people, see, for example, Collins [Co], He [He], He and Sun [HS1][HS2], and Huang [Hu].

Recently Bedulli, He and Vezzoni [BHV] gave an alternative approach  to the existence (resp. uniqueness) of the flow (1.1) using the transverse De Turck trick (resp. Kotschwar's idea). They also considered the transverse K$\ddot{a}$hler-Ricci flow for  K$\ddot{a}$hler foliations, which is both a special case of (1.1) and a generalization of (1.2).

Now let $\mathcal{F}$ be a  transversely holomorphic  foliation of real codimension $2n$  on a compact manifold $M^{m+2n}$. Suppose $\mathcal{F}$ is generated by the involutive distribution $L$ of rank $m$. Let $Q=TM/L$ be the normal bundle to the foliation. We have the following exact sequence:
\begin{equation}
0\rightarrow L \rightarrow TM \rightarrow  Q \rightarrow 0.
\end{equation}
A transversely Hermitian structure on $\mathcal{F}$ is given by the transverse complex structure $J$ and a Hermitian metric $g^T$  ($g^T(J\cdot,J\cdot)=g^T(\cdot,\cdot)$) on the  bundle $Q$ satisfying $\mathcal{L}_Xg^T=0$ for any $X\in \Gamma (L)$.   Let $\omega$ be the pull-back of $g^T(J\cdot, \cdot)$ to $M$. For convenience we'll call $\omega$ a transversely Hermitian metric.  We say  such $\omega$ is positive, denoted by $\omega>0$.  In local complex foliated coordinates  $(x^1,\cdot\cdot\cdot,x^m, z^1, \cdot\cdot\cdot,z^n)$  (see [BHV]), we may write  $g^T=g^T_{i\bar{j}}dz^id\bar{z}^j$. Then  the expression
$-\sqrt{-1}\partial_B \bar{\partial}_B \log \det (g^T_{i\bar{j}})$ does not depend on the local complex foliated coordinates, and gives a closed basic 2-form on $M$, denoted by $\rho^T$.  We call $\rho^T$ the transverse Chern-Ricci form. Its basic cohomology class in the basic Bott-Chern cohomology group
\begin{equation*}
H^{1,1}_{BBC}(M,\mathbb{R})=\frac{\{\text {closed basic real (1,1)-forms}\}}{\{\sqrt{-1}\partial_B \bar{\partial}_B \psi, \psi \in C^\infty_B(M,\mathbb{R})\}}
\end{equation*}
is the basic first Bott-Chern class, denoted by $c_1^{BBC}(M)$.

 In this note we'll consider the following transverse Chern-Ricci flow for transversely Hermitian foliations on compact manifolds,
\begin{equation}
\frac{\partial  \omega(t)}{\partial t}=-\rho^T(\omega(t)), \hspace*{8mm}
\omega(0)=\omega_0,
\end{equation}
which is a generalization of the transverse K$\ddot{a}$hler-Ricci flow, and is analogous to the Chern-Ricci flow ([G][TW2]). Here $\omega(t)$ is a family of transversely Hermitian metrics, and $\rho^T(\omega(t))$ is the transverse Chern-Ricci form of $\omega(t)$.

\begin{thm} \label{thm 1.1} \ \
  Let $\mathcal{F}$ be a transversely holomorphic foliation on a compact manifold $M$ with a transversely Hermitian metric $\omega_0$.  Then the flow (1.4) has a unique maximal time solution on $[0, T_0)$, where
  $T_0=\sup \{t\geq 0| \hspace*{1mm} \exists \psi \in C^\infty_B(M) \hspace*{1mm} with \hspace*{1mm} \omega_0-t\rho^T(\omega_0)
  +\sqrt{-1}\partial_B \bar{\partial}_B \psi >0\}$.
\end{thm}

For the K$\ddot{a}$hler case of Theorem 1.1 see Tian-Zhang [TZ], the Hermitian case see  Tosatti-Weinkove [TW2], and the transverse   K$\ddot{a}$hler case see   Bedulli-He-Vezzoni [BHV].

The following theorem is analogous to Theorem 1.7 in [TW2], and extends Theorem 6.3 in [BHV].

\begin{thm} \label{thm 1.2} \ \
  Let $\mathcal{F}$ be a transversely holomorphic foliation on a compact manifold $M$ with a transversely Hermitian metric $\omega_0$.  Suppose that  $c_1^{BBC}(M) <0$. Then the flow (1.4)  exists for all time and as $t\rightarrow \infty$  the rescaled transversely Hermitian
metrics $\omega(t)/t$ converge smoothly to a transverse K$\ddot{a}$hler-Einstein metric. 
\end{thm}

Here the condition $c_1^{BBC}(M) <0$ means that the class $-c_1^{BBC}(M)$ is represented by a positive, closed basic real $(1,1)$-form, hence a transverse K$\ddot{a}$hler form.

Recall that a transversely oriented Riemannian foliation $\mathcal{F}$ is homologically orientable if there is a $m$-form  $\chi$ on $M$ which restricts to a volume form on each leaf and satisfies $d\chi(X_1,\cdot\cdot\cdot,X_m,Y)=0$, where $X_1,\cdot\cdot\cdot,X_m$ are tangent to the leaves; see for example [EK2].

\begin{thm} \label{thm 1.3} \ \
   Let $\mathcal{F}$ be a transversely holomorphic foliation on a compact manifold $M$ with a transversely Hermitian metric $\omega_0$. Suppose that $\mathcal{F}$ is homologically orientable and that the basic first Bott-Chern class of $(M,\mathcal{F})$ is zero. Then the flow (1.4) exists for all time and as $t\rightarrow \infty$  converges smoothly to a transversely Hermitian metric $\omega_\infty$ with $\rho^T(\omega_\infty)=0$.
\end{thm}

For the K$\ddot{a}$hler case of Theorem 1.3 see Cao [C],  the Hermitian case see Gill [G],  the transverse K$\ddot{a}$hler case see Bedulli-He-Vezzoni [BHV].

The short time existence of the solution to the transverse Chern-Ricci flow follows from a general existence result on  strongly transversely   parabolic systems of arbitrary order, see Section 2.
In Section 3 we prove Theorems 1.1, 1.2 and 1.3, following [TW2] and [G].

\section{ Strongly transversely parabolic systems}

Let $\mathcal{F}$ be a transversely oriented Riemannian foliation of codimension $q$ on a compact manifold $M$, and
$(E, \nabla)$ be a  $\mathcal{F}$-vector bundle over $M$. (See for example  [EK2].)
Let
\begin{equation}
D: C^\infty(E/\mathcal{F}) \rightarrow C^\infty(E/\mathcal{F})
\end{equation}
be a  (smooth) basic differential operator  of order $r$ of the form $ D(u)=F(x, u, \nabla u, \cdot\cdot\cdot, \nabla^{r}u)$, where $x\in M$. ($F$ may be viewed as a map from  $C^\infty (J^r(E/\mathcal{F}))$ to $C^\infty(E/\mathcal{F})$; compare Lemma 2.1 below.) Here `basic' means that in local foliated coordinates $(x^1, \cdot\cdot\cdot, x^m, y^1, \cdot\cdot\cdot, y^q)$ and local frame $\{e_a\}_{a=1}^l$  of $E$,
\begin{equation}
  D(u)=F^a(y^1,\cdot\cdot\cdot,y^q, u^1,\cdot\cdot\cdot,u^l, \frac{\partial u^1}{\partial y^1},\cdot\cdot\cdot, (\frac{\partial}{\partial y^q})^{r} u^l)e_a,
\end{equation}
where $F^a$ are smooth functions of their arguments.

Given $w, v\in C^\infty(E/\mathcal{F})$, the linearization of the operator $D$ at $w$ in the direction $v$ is
$ D_{*|w}(v)=\frac{\partial}{\partial s}(D(w+sv))|_{s=0}=\lim_{s\rightarrow  0}\frac{D(w+sv)-D(w)}{s}$. Given  a transverse covector $(x, \xi)\in T^*M$ (that is, $\xi(X)=0$ for any $X\in T_xM$ tangent to the leaf), choose $\phi \in C^\infty_B(M)$ with $d\phi(x)=\xi$.  Define the principal symbol $\sigma(D_{*|w})$ via
\begin{equation}
\sigma(D_{*|w})(x,\xi)v=\lim_{s\rightarrow \infty}s^{-r}e^{-s\phi(x)}D_{*|w}(e^{s\phi}v)(x).
\end{equation}
Compare Topping [T] and [EK1]. The operator $D$ is strongly transversely elliptic at $w\in C^\infty(E/\mathcal{F})$  if $r$ is even and there exists a constant $\mu>0$ such that
\begin{equation}
(-1)^{r/2} \langle\sigma(D_{*|w})(x,\xi)v, v\rangle \geq \mu |\xi|^{r}|v|^2
\end{equation}
for all transverse covector $(x, \xi)\in T^*M$ and $v\in C^\infty(E/\mathcal{F})$; here   $\langle\cdot,\cdot\rangle$ is some fiber metric on $E$.

Now, following [EK1] (see also [BHV]) we can form the basic $r$-jet bundle $J^r(E/\mathcal{F})$, whose fiber over $x\in M$ is $ C^\infty (E/\mathcal{F})/Z_x^r(E/\mathcal{F})$,  where $Z_x^r(E/\mathcal{F})$ is the ring of basic sections $u$ of $E$ satisfying $(\nabla^k u)(x)=0$ for all $0\leq k \leq r$, and there is a natural map $J_r: C^\infty (E/\mathcal{F}) \rightarrow C^\infty (J^r(E/\mathcal{F}))$.

\begin{lem} \label{lem 1.2} \ \
Given a basic differential operator $D$ of order $r$ as above, there exists  a foliated map $T: J^r(E/\mathcal{F})\rightarrow E$ such that $D=T_* \circ J_r$, where $T_*: C^\infty (J^r(E/\mathcal{F})) \rightarrow C^\infty (E/\mathcal{F})$ is the map induced by $T$.
\end{lem}
{\bf Proof}.   It is a slight adaption of proof of Theorem 1 on pp. 61-62 in [P].  By definition $D$ maps $Z_x^r(E/\mathcal{F})$ to a singe point in $E_x$.  So there is a unique map $T(x): C^\infty (E/\mathcal{F})/Z_x^r(E/\mathcal{F}) \rightarrow E_x$ such that $D(f)(x)=T(x)J_r(f)_x$. Since $D$ is basic, $T$ is foliated. Clearly $D=T_* \circ J_r$.   \hfill{$\Box$}

\vspace*{0.4cm}

Let $G=SO(q)$ and $\rho: M^{\sharp} \rightarrow M$ be the $G$-principal bundle of oriented  orthonormal frames transverse to $\mathcal{F}$. Then, by Molino [M], the foliation  $\mathcal{F}$ can be lifted to a $G$-invariant foliation $\mathcal{F}^{\sharp}$ on $M^{\sharp}$ which is transversely parallelizable; furthermore, there are a manifold $W$ and a locally trivial fibration
$\pi^{\sharp}: M^{\sharp} \rightarrow W$, called the basic fibration of $\mathcal{F}$,  whose fibres are the closures of the leaves of $\mathcal{F}^{\sharp}$, which are submanifolds of $M^{\sharp}$.

  We pull-back $E$ via $\rho$ to a bundle $E^{\sharp}$ over $M^{\sharp}$. For $u \in W$, let $F^{\sharp}_u$ be the fiber over $u$ of the basic fibration $\pi^{\sharp}$, and $E^{\sharp}_u$ (resp. $\mathcal{F}^{\sharp}_u$) be the restriction of $E^{\sharp}$ (resp. $\mathcal{F}^{\sharp}$) to $F^{\sharp}_u$.
Let $\bar{E}_u=C^\infty (E^{\sharp}_u/\mathcal{F}^{\sharp}_u)$, we can paste $\bar{E}_u$ together to form a Hermitian vector bundle $\bar{E}$ over $W$, such that there are canonical isomorphisms
$\psi^{\sharp}: C^\infty(E^{\sharp}/\mathcal{F}^{\sharp}) \rightarrow C^\infty(\bar{E})$ and $\psi: C^\infty(E/\mathcal{F}) \rightarrow C_G^\infty(\bar{E})$; see [EK1].

The following theorem extends [EK1] and [BHV], which consider the linear  and quasilinear cases respectively.

\begin{lem} \label{lem 2.2} \ \ Let $D: C^\infty (E/\mathcal{F}) \rightarrow C^\infty (E/\mathcal{F})$ be a basic differential operator which is strongly transversely elliptic at $u\in C^\infty (E/\mathcal{F})$. Then there is a  differential operator $\bar{D}: C^\infty_G(\bar{E}) \rightarrow C^\infty_G(\bar{E})$ which is strongly transversely elliptic at $\psi(u)$ and satisfies $D\circ \psi=\psi\circ \bar{D}$.
\end{lem}

{\bf Proof}   The proof is almost the same as that in [EK1] and [BHV].   Let $D: C^\infty (E/\mathcal{F}) \rightarrow C^\infty (E/\mathcal{F})$ be a basic differential operator of order $r$. Following [EK1] (see also [BHV]),  from the map $T$ obtained  in Lemma 2.1 we can construct a map $T^{\sharp}: \bigoplus_{k=0}^r S^k(Q^{\sharp}, E^{\sharp}) \rightarrow E^{\sharp}$ such that
$\rho^{\sharp} \circ T^{\sharp}=T\circ \rho^{\sharp} _s: \bigoplus_{k=0}^r S^k(Q^{\sharp}, E^{\sharp}) \rightarrow E$, where $\rho^{\sharp} _s: \bigoplus_{k=0}^r S^k(Q^{\sharp}, E^{\sharp}) \rightarrow \bigoplus_{k=0}^r S^k(Q, E) \cong J^r(E/\mathcal{F})$ is the map induced by $\rho^{\sharp}$.

Let $D^{\sharp}=T^{\sharp}_*\circ J_r^{\sharp}: C^\infty (E^{\sharp}/\mathcal{F}^{\sharp}) \rightarrow C^\infty (E^{\sharp}/\mathcal{F}^{\sharp}) $.

Suppose $D$ is strongly transversally elliptic at $u$.

 As in [BHV],  differentiating the equality $\rho^{\sharp} \circ T^{\sharp}=T\circ \rho^{\sharp} _s: C^\infty (J^r(E^{\sharp}/\mathcal{F}^{\sharp})) \rightarrow C^\infty (E)$ at $u^{\sharp}$,  we get $\rho^{\sharp} \circ T^{\sharp}_{*|J_r(u^{\sharp})}=T_{*|J_r(u)}\circ \rho^{\sharp} _s: C^\infty (J^r(E^{\sharp}/\mathcal{F}^{\sharp})) \rightarrow C^\infty (E)$.

As in [EK1], let $D'=D^{\sharp}+(-1)^{r/2}(\sum_{j=1}^N Q_j\circ Q_j)^{r/2}$, where $N=\frac{1}{2}q(q-1)$ and $Q_j:C^\infty (E^{\sharp}) \rightarrow C^\infty (E^{\sharp})$ are induced from a basis of Lie algebra of $G=SO(q)$.
Finally let $\bar{D}=\psi^{\sharp}D'(\psi^{\sharp})^{-1}: C^\infty(\bar{E}) \rightarrow  C^\infty(\bar{E}) $. By construction
$\bar{D}$ is $G$-invariant, strongly transversally elliptic at $\psi(u)$, and  satisfies $\psi \circ D=\bar{D} \circ \psi: C^\infty (E/\mathcal{F}) \rightarrow C_G^\infty(\bar{E})$.
\hfill{$\Box$}

\vspace*{0.4cm}

Let $u_0$  be a given smooth basic function on $M$, consider the equation
\begin{equation}
\frac{\partial u}{\partial t}+D(u(\cdot,t))=0,  \hspace*{2mm} u(\cdot,0)=u_0,
\end{equation}
for $u(\cdot,t) \in C^\infty_B(M)$.

As preparation we first treat  the non-foliated case, that is, the case $m=0$; compare the Main Theorem 1  in Baker [B].
\begin{thm} \label{thm 2.2} \ \ Let $m=0$.
Suppose $D$ is strongly  elliptic at $u_0$. Then the equation (2.5) has a unique short time solution.
\end{thm}
{\bf Proof}.   We follow the proof of Main Theorem 1 in Baker [B], cf. also Hamilton [H1].  By the Schauder theory for the linear, strongly parabolic systems (see for example Friedman [F1] [F2],  Eidelman [E], and Lamm [L]), the linearization of (2.5) at $u_0$,
\begin{equation}
\frac{\partial w}{\partial t}+D_{*|u_0}(w(\cdot,t))=0,  \hspace*{2mm} w(\cdot,0)=u_0,
\end{equation}
is uniquely solvable. Then we apply the Schauder theory to the linearization of (2.5) at the solution $w(\cdot,t)$ to (2.6),
\begin{equation}
\frac{\partial v}{\partial t}+D_{*|w(\cdot,t)}(v(\cdot,t))=0,  \hspace*{2mm} v(\cdot,0)=u_0,
 \end{equation}
 and see that the Fr$\acute{e}$chet derivative of the operator $\frac{\partial}{\partial t}+D$ is invertible at $w$. Finally by the Schauder estimate and inverse function theorem we can find a unique solution to (2.5) for a short time  near the solution $w$ to (2.6).  See [B] for more details. \hfill{$\Box$}

\vspace*{0.4cm}

{\bf Remark} The condition 1) in Main Theorem 1 of [B] is not essential, since it is not needed in the Schauder theory for the linear, strongly parabolic systems; see for example Friedman [F1] [F2],  Eidelman [E], Lamm [L]   and Section 10 of Chapter VII in [LSU]; compare also [GM].
The conditions 3) and 4) in Main Theorem 1 of [B]  are automatically satisfied in our situation, since the differential operator $D$ that we consider is smooth, and the manifold $M$ is compact.

\vspace*{0.4cm}

Now we treat the foliated case. For the quasilinear case, see Theorem 1.1 in [BHV].
\begin{thm} \label{thm 2.3} \ \  Suppose $D$ is strongly transversely elliptic at $u_0$. Then the equation (2.5) has a unique short time solution.
\end{thm}

{\bf Proof}. With the help of Lemma 2.2 and Theorem 2.3, the proof is similar to that of Theorem 1.1 in [BHV]. \hfill{$\Box$}

\section{ Proof of Theorems 1.1, 1.2 and 1.3}

\vspace*{0.4cm}

\begin{lem} \label{lem 3.1} \ \ Let $\mathcal{F}$ be a transversely Hermitian foliation on a compact manifold $M$.  There exists a unique connection on the normal bundle $Q$ which is adapted to the Bott connection (see for example [LMR]) and which is compatible with the transverse complex structure and the transverse Hermitian metric.
\end{lem}
{\bf Proof}.  Consider  local foliated charts $\{U_i\}$ which constitute a locally finite, countable cover of $M$, and  submersions $f_i: U_i \rightarrow T$, where $T$ is the  transverse manifold of the foliation (see [BE]).  Now `pull-back' the Chern connection on the Hermitian manifold $T$ via $f_i$ and patch them together using partition of unity subordinate to the cover $\{U_i\}$.  \hfill{$\Box$}

\vspace*{0.4cm}

We may call the connection in Lemma 3.1 transverse Chern connection for the Hermitian foliation $\mathcal{F}$. (There are other ways to prove Lemma 3.1.)

\vspace*{0.4cm}

Observe that Proposition 6.9 (maximum principle for basic maps) and its proof in [BHV] hold true in the more general situation of transversely Hermitian foliations.

\vspace*{0.4cm}

With the help of Lemmas 3.1 and the above observation,  we can  prove Theorem 1.1 by adapting the proof of Theorem 1.2 in [TW2].
Let $\alpha_t=\omega_0-t\rho^T(\omega_0)$.
  Fix $T'< T_0$.  Then there is a basic function $f_{T'}$ with $\alpha_{T'}+\sqrt{-1}\partial_B \bar{\partial}_Bf_{T'}>0$.  Let
\begin{equation}
\lambda=\frac{1}{T'}\sqrt{-1}\partial_B \bar{\partial}_Bf_{T'}-\rho^T(\omega_0),
\end{equation}
and
 \begin{equation}
\hat{\omega}_t=\frac{T'-t}{T'}\omega_0+\frac{t}{T'}(\alpha_{T'}+\sqrt{-1}\partial_B \bar{\partial}_Bf_{T'}).
\end{equation}
Then  $\hat{\omega}_t=\omega_0+t\lambda$.
Choose a transverse volume form $\Omega:=e^{\frac{f_{T'}}{T'}}\omega_0^n$. Then $\sqrt{-1}\partial_B \bar{\partial}_B \log \Omega=\lambda$.
Consider the transversely parabolic complex Monge-Amp$\grave{e}$re equation
\begin{equation}
\frac{\partial}{\partial t}\varphi= \log \frac{(\hat{\omega}_t+\sqrt{-1}\partial_B \bar{\partial}_B\varphi)^n}{\Omega},  \hspace*{2mm} \hat{\omega}_t+\sqrt{-1}\partial_B \bar{\partial}_B\varphi>0,  \hspace*{2mm}  \varphi|_{t=0}=0
\end{equation}
for basic functions $\varphi(\cdot, t)$.

Given a solution $\varphi(\cdot,t) \in C^\infty_B(M)$ to (3.3), we see that $\omega_t:=\hat{\omega}_t+\sqrt{-1}\partial_B \bar{\partial}_B\varphi$ is a solution to
(1.4) by taking $\sqrt{-1}\partial_B \bar{\partial}_B$ of (3.3). On the other hand, given a solution $\omega_t$ to (1.4) on some subinterval of $[0,T']$, one computes
\begin{equation}
\frac{\partial}{\partial t}(\omega_t-\hat{\omega}_t)=\sqrt{-1}\partial_B \bar{\partial}_B\log \frac{\omega_t^n}{\Omega}.
\end{equation}
 Then for any fixed point $x$ on $M$ we solve the ODE
 \begin{equation}
 \frac{d}{dt}\varphi(x,\cdot)=\log \frac{\omega(x,\cdot)^n}{\Omega(x)}, \hspace*{2mm} \varphi(x,0)=0.
 \end{equation}
 Clearly $\varphi(\cdot,t)$ is basic. From (3.4) and (3.5) it is easy to see that $\omega_t=\hat{\omega}_t+\sqrt{-1}\partial_B \bar{\partial}_B\varphi$, and $\varphi$ solves (3.3).

 Clearly, the RHS of (3.3) is strongly transversely elliptic  at $\varphi=0$. So by Theorem 2.4, (3.3) has a unique short time solution. (Actually here we only need the (fully nonlinear) second order equation case of Theorem 2.4, whose proof is much easier. In this case we only need the Schauder theory for the linear parabolic equation, see for example [LSU] and [Li].  Anyway Theorem 2.3 in the  second order equation case is classical.)

Now the rest of the arguments is almost the same as that in [TW2].  Note that instead of developing a transverse Evans-Krylov theory we may easily adapt the third  order estimate in Cherrier [Ch] to the foliated case.  \hfill{$\Box$

\vspace*{0.4cm}

The proof of Theorem 1.2 is almost the same as that of Theorem 1.7 in [TW2].  In the proof we use the normalized transverse Chern-Ricci flow,
\begin{equation*}
\frac{\partial  \omega(t)}{\partial t}=-\rho^T(\omega(t))-\omega(t), \hspace*{4mm}
\omega(0)=\omega_0,
\end{equation*}
and the foliated maximum principle.
\hfill{$\Box$}

\vspace*{0.4cm}

Now we assume that $\mathcal{F}$ is a transversely holomorphic foliation (of real codimension $2n$) on a compact manifold $M$ and is homologically orientable, so  there is a $m$-form  $\chi$ on $M$ which restricts to a volume form on each leaf and which is $\mathcal{F}$-relatively closed (see Introduction).

We need a lemma on integration by parts, which is well-known.
\begin{lem} \label{lem 3.2} \ \ Let $\mathcal{F}$  and $\chi$ be as above.
Let $\alpha$ and $\beta$ be basic forms with $\deg \alpha +\deg \beta=2n-1$. Then
\begin{equation}
\int_M d_B\alpha \wedge \beta \wedge \chi=(-1)^{\deg \alpha +1}\int_M\alpha \wedge d_B\beta \wedge \chi.
 \end{equation}
\end{lem}
{\bf Proof}. Since $\alpha$ and $\beta$ are basic forms and  $\deg \alpha +\deg \beta=2n-1$, using the definition of $\chi$ we see that $\alpha \wedge \beta \wedge d\chi=0$. Then
\begin{equation}
d(\alpha \wedge \beta \wedge \chi) = d_B\alpha \wedge \beta \wedge \chi+(-1)^{\deg \alpha}\alpha \wedge d_B\beta  \wedge \chi,
\end{equation}
and the lemma follows. \hfill{$\Box$}

\vspace*{0.4cm}
\begin{lem} \label{lem 3.3} \ \  Let $\mathcal{F}$ be a transversely holomorphic foliation (of real codimension $2n$) on a compact manifold $M$ with a transversely Hermitian metric $\omega_0$. Suppose that $\mathcal{F}$ is homologically orientable. Then there is a  basic function $u$  such that the transverse Hermitian metric $\omega_G=e^u\omega_0$ is transversely Gauduchon, i.e., $\partial_B \bar{\partial}_B (\omega_G^{n-1})=0$.
\end{lem}
{\bf Proof}.  The proof of Theorem (1.2.4) (Gauduchon's theorem) on p. 30 in [LT] is easily adapted to our situation. 
\hfill{$\Box$}

\vspace*{0.4cm}

Now we can adapt Gill [G] to prove Theorem 1.3.  The proof is very similar.  We only point out that for the zeroth order estimate we use Lemma 3.2, Lemma 3.3, and the Moser iteration; compare [TW1] and [SWZ]. \hfill{$\Box$}

\vspace*{0.4cm}

\noindent {\bf Acknowledgements} {\hspace*{4mm}}   I'm partially supported by NSFC no.11171025.


\hspace *{0.4cm}

\bibliographystyle{amsplain}

\noindent {\bf Reference}

\hspace *{0.1cm}

\vspace*{0.4cm}

[B] C. Baker, The mean curvature flow of submanifolds of high codimension,  arXiv:1104.4409.

[BE] R. Barre,  A. El Kacimi-Alaoui,  Foliations, in  Handbook of differential geometry, Vol. II,  35-77, Elsevier/North-Holland, Amsterdam, 2006.

[BHV] L. Bedulli, W. He, L. Vezzoni, Second order geometric flows on foliated manifolds,   arXiv:1505.03258.

[C] H.D. Cao, Deformation of K$\ddot{a}$hler metrics to K$\ddot{a}$hler-Einstein metrics,  on compact K$\ddot{a}$hler manifolds, Invent. Math. 81 (1985), no. 2, 359-372.

[Ch] P. Cherrier, $\acute{E}$quations de Monge-Amp$\grave{e}$re sur les vari$\acute{e}$t$\acute{e}$s Hermitiennes compactes, Bull. Sc. Math. (2) 111 (1987), 343-385.

[Co] T. Collins, The transverse entropy functional
and the Sasaki-Ricci flow, Trans. Amer. Math. Soc. 365 (2013),
1277-1303.

[E] S. Eidelman, Parabolic systems,  North-Holland Publishing Co., 1969.

[EK1] A. El Kacimi-Alaoui, Op¨¦rateurs transversalement elliptiques sur un feuilletage riemannien et applications, Compositio Math. 73 (1990), no. 1, 57-106.

[EK2] A. El Kacimi-Alaoui, Towards a basic index theory, in  Dirac operators: yesterday and today,  251-261, Int. Press, 2005.

[F1] A. Friedman,  Interior estimates for parabolic systems of partial differential equations, J. Math. Mech. 7 (1958), 393-417.

[F2] A. Friedman, Partial differential equations of parabolic type, Prentice-Hall, Inc., 1964.

[GM] M. Giaquinta, G. Modica, Local existence for quasilinear parabolic systems under nonlinear boundary conditions, Ann. Mat. Pura Appl. 4 149 (1987), 41-59.

[G] M. Gill, Convergence of the parabolic complex Monge-Amp$\grave{e}$re equation on compact Hermitian manifolds, Comm. Anal. Geom. 19 (2011), no. 2, 277-303.

[H1] R. Hamilton, Harmonic maps of manifolds with boundary, LNM 471, Springer, 1975.

[H2] R. Hamilton, Three-manifolds with positive Ricci curvature, J. Diff. Geom. 17 (1982), no. 2, 255-306.

[He] W. He, The Sasaki-Ricci flow and compact
Sasaki manifolds of positive transverse holomorphic bisectional
curvature, J. Geom. Anal. 23 (2013), 1876-1931.

[HS1] W. He, S. Sun,  Frankel conjecture and Sasaki Geometry, arXiv:1202.2589.

[HS2] W. He, S. Sun, The generalized Frankel conjecture in Sasaki
geometry, Inter. Math. Res. Not. 2013,  doi:10.1093/imrn/rnt185

[Hu]  H. Huang, Sasaki manifolds with positive transverse orthogonal bisectional curvature, arXiv:1301.1229, to appear in Advances in Geometry (2015).

[LSU] O. A. Lady$\check{z}$enskaja, V. A. Solonnikov, N. N. Ural'ceva, Linear and quasilinear equations of parabolic type, Amer. Math. Soc., 1968.

[L] T. Lamm, Biharmonischer W$\ddot{a}$rmefluss, Diploma thesis,  Albert-Ludwigs-Universit$\ddot{a}$t Freiburg, 2002.

[Li]  G. Lieberman, Second order parabolic differential equations,  World Scientific Publishing Co., Inc., 1996

[LMR] M. Lovri$\acute{c}$, M. Min-Oo, E. Ruh, Deforming transverse Riemannian metrics of foliations, Asian J. Math. 4 (2000), no.2, 303-314.

[LT] M. L$\ddot{u}$bke, A. Teleman, The Kobayashi-Hitchin correspondence, World Scientific, 1995.

[M] P. Molino, Riemannian foliations,  Birkh$\ddot{a}$user, 1988.

[P] R. Palais, Seminar on the Atiyah-Singer index theorem, Ann. Math. Studies no. 57, Princeton University Press, 1965.

 [TZ] G. Tian, Z. Zhang, On the K$\ddot{a}$hler-Ricci flow on projective manifolds of general type, Chinese Ann. Math. Ser. B 37 (2006), no. 2, 179-192.

[T] P. Topping, Lectures on the Ricci flow, Cambridge University Press, 2006.

[TW1] V. Tosatti, B. Weinkove, The complex Monge-Amp$\grave{e}$re equation on compact Hermitian manifolds, J. Amer. Math. Soc. 23 (2010), no. 4, 1187-1195.

[TW2] V. Tosatti, B. Weinkove,  On the evolution of a Hermitian metric by its Chern-Ricci form, J. Diff. Geom. 99 (2015), 125-163.

[SWZ] K. Smoczyk, G. Wang, Y. Zhang, The Sasaki-Ricci flow, Internat. J.
Math. 21 (2010), no. 7, 951-969.

\vspace *{0.4cm}

School of Mathematical Sciences, Beijing Normal University,

Laboratory of Mathematics and Complex Systems, Ministry of Education,

Beijing 100875, P.R. China

 E-mail address: hhuang@bnu.edu.cn

\end{document}